\documentclass[3p,times]{elsarticle}

\usepackage[cp1251]{inputenc}%�������� ���� ���������:  utf8 � UNIX, cp1251 � Windows

\usepackage[english]{babel}
\usepackage{amsthm,amssymb,amsmath,amsfonts}
\usepackage{graphicx}
\usepackage{hhline,array}
\usepackage{color}
\usepackage{tabularx}

\newtheorem{theorem}{\bf Theorem}[section]

\usepackage{ecrc}

\volume{00}

\firstpage{1}

\journalname{ }

\runauth{A. N. Tynda, D. N. Sidorov, N. A.Sidorov}

\jid{ }

\usepackage{amssymb}
\usepackage[figuresright]{rotating}

\begin{document}

\begin{frontmatter}



\dochead{}

\title{Numeric solution of systems of nonlinear Volterra integral equations of the first kind with discontinuous kernels}

\author[label1]{Aleksandr Tynda}
\author[label2]{Denis Sidorov\corref{corauth}}
\cortext[corauth]{Corresponding author.}
\ead{contact.dns@gmail.com,  {telephone number: +79148822035}}
\author[label3]{Nikolai  Sidorov}
%DS was funded NSFC Grant No. 61673398.}

\address[label1]{Penza State University}
\address[label2]{Energy Systems Institute of Russian Academy of Sciences}
\address[label3]{Irkutsk State University}

\begin{abstract}

The systems of nonlinear Volterra integral equations of the first kind with jump discontinuous kernels are studied.  The iterative numerical method for such nonlinear systems
is proposed. Proposed method employs the modified Newton-Kantorovich iterative process  for the integral operators linearization. On each step of the iterative
process the linear system of integral equations is obtained and resolved using the
discontinuity driven piecewise constant  approximation of the exact solution. The convergence theorem
is proved.  The polynomial collocation technique is also applied to solve such systems  as the alternative method. The accuracy of proposed numerical methods is discussed. The  model examples are examined in order to demonstrate the efficiency of proposed numerical methods and illustrate the constructed theory.

\end{abstract}


\end{frontmatter}

% \linenumbers

%% mai

%##########################################################################################
%******************************************************************************************
%---------------------------     S E C T I O N  1 -----------------------------------------
%******************************************************************************************
%##########################################################################################

\section*{Introduction}

Functional equations such as differential, integral and integral-differential with various delay types are univer\-sal tools for dynamical systems modeling in
physics, engineering, medical and ecological systems, economics and other fields, see  \cite{kantor}, \cite{hriton}, \cite{apartsyn},  \cite{Sidorov-book} and their
bibliography. Models based on such mathematical tools are enable to properly describe the process of interest and its the only way to do it properly
in most of the cases. Indeed, integral equations with delays able to describe various dynamical processes. One of the most interesting classes of functional equations with delays
is the Volterra (or evolutionary) integral equations (VIE). Volterra models take into account the memory effects of the system when the past states influence the present.
Integral equations are useful tool for dynamical models modeling with various applications \cite{arxiv}. It is to be noted that exact solutions of such equations can not be derived in most of the practical problems and it is necessary to develop the efficient numerical methods.

The Volterra equation is classical problem have been considered by numerous studies \cite{book}. But only few authors have considered these
equations with discontinuous kernels.
G.C.~Evans \cite{Evans} first considered the
integral {e}quation of the {s}econd {k}ind with discontinuous {k}ernel.
P.M.~Anselone \cite{Anselone} proposed the uniform approximation {t}heory for integral {e}quations with
{d}iscontinuous {k}ernels. A.M.~Denisov and A.~Lorenzi  \cite{denisovlorenzi} derived the sufficient conditions of solution existence of the linear VIE
with jump discontinuity on the single curve. The excellent overview of the development of theory of the Volterra integral equations
of the first kind during the last century is given by H.~Brunner \cite{brunner}.

The theory of the linear VIE of the 1st kind
 with discontinuous kernels along the finite number of continuous curves was first considered in \cite{Sidorov-FR, SidorovDU}, pursued in articles \cite{Sidorov-2013,SidorovMarkovaAiT,vestnik} followed by monograph \cite{Sidorov-book}.
 A.~Lorenzi \cite{Lorenzi2013} and N.A.~Sidorov \& D.N.~Sidorov \cite{sidsid} constructed the theory of VIE with discontinuous kernels in Banach spaces applicable to the wide classes of integral-differential PDEs.
 A.P.~Khromov \cite{Khromov} considered  the integral {o}perators with {d}iscontinuous {k}ernel on {p}iecewise {l}inear {c}urves. P.J.~Davies and D.B.~Duncan \cite{Davies} employed the cubic ``convolution spline'' method for
numerical approximation of first kind Volterra convolution integral equations with discontinuous kernels.
%J. Integral Equations Appl. 29, No. 1, 41-73 (2017).
 In the series of articles \cite{Boikov-Tynda11} -- \cite{Muftahov-CAM} there were number of
efficient numerical methods have been proposed for solution of linear and nonlinear equations
with variable limits of integration which can be interpreted as time delays.

 In this paper we generalize our results and construct the efficient numerical  methods of solution to nonlinear systems of VIE with discontinuous kernels based on direct quadratures.

%***********to be edited!****************
The paper is organized into five sections. In Section 1  the
problem statement is formulated.
Section 2 is dedicated to the modified Newton-Kantorovich iterative process' application to the nonlinear systems of VIEs . Section 3 deals with the discretization process of the system of linear integral equation appearing on each step of
proposed iterative process in Section 2.
The constructed theory and numerical methods are illustrated on model examples in Section 4.
The concluding remarks  are given in Section 5.
%***************************

%##########################################################################################
%******************************************************************************************
%---------------------------     S E C T I O N  2 -----------------------------------------
%******************************************************************************************
%##########################################################################################

\section{Problem statement}

Let us consider the following nonlinear system of the Volterra integral equations
of the first kind
\begin{equation}
\label{MainSystem}
  \begin{cases}
     \int\limits_{0}^{t} h_{1}(t,s,x_1(s),x_2(s), ..., x_n(s))ds - f_1(t) = 0,\\
     \int\limits_{0}^{t} h_{2}(t,s,x_1(s),x_2(s), ..., x_n(s))ds - f_2(t) = 0,\\
     \vdots\\
     \int\limits_{0}^{t} h_{n}(t,s,x_1(s),x_2(s), ..., x_n(s))ds - f_n(t) = 0.
\end{cases}
\end{equation}
where \(t\in [0,T]\) and kernels \(h_{i}(t,s,x_1(s),x_2(s), ..., x_n(s))\) are jump discontinuous
on the curves
\(\alpha_j(t),\) \(j = \overline{1,n-1},\)  and can be defined as
\begin{equation}\label{Kernels}
   h_{i}(t,s,x_1(s),x_2(s), ..., x_n(s)) =
   \begin{cases}
     K_{i1}(t,s)G_{i1}(s,x_1(s)), \; (t,s) \in m_1,\\
     K_{i2}(t,s)G_{i2}(s,x_2(s)), \; (t,s) \in m_2,\\
     \vdots\\
     K_{in}(t,s)G_{in}(s,x_n(s)), \;  (t,s) \in m_n.
   \end{cases}
\end{equation}
Here \(m_j = \{(t,s)|\;\alpha_{j-1}(t) < s \le \alpha_j(t)\}\), \(\alpha_0(t) = 0\), \(\alpha_n(t) = t\), \(i = \overline{1,n}\), \(j = \overline{1,n}\). Functions
\(K_{ij}(t,s)\), \(f_i(t)\), \(\alpha_j(t)\) enjoy continuous derivatives wrt time \(t\) for \((t,s) \in \overline{m}_j\), \(K_{in}(t,t) \ne 0\),  \(0 < \alpha_1(t) < \alpha_2(t) < ... < \alpha_{n-1}(t) < t\). Functions \(\alpha_1(t), ..., \alpha_{n-1}(t)\) are nondecreazing  for \(t\in [0,T]\).

For the statements concerning the existence and uniqueness of solutions of systems of VIEs  with discontinuous kernels readers may refer to \cite{Muftahov-Ural}.

In this paper we propose an iterative method for solving such nonlinear systems based on linearization of the integral vector-operator. To solve the systems of linear equations with discontinuous kernels arising in the iterative process, we suggest two approaches: a generalization of the direct discretization algorithm proposed in \cite{Muftahov-CAM} and the polynomial-collocation scheme. The first approach is based on piecewise constant approximation of exact solutions and has the first order of accuracy. The second approach uses the polynomial approximation of exact solutions and allows to derive more accurate solutions.

%##########################################################################################
%******************************************************************************************
%---------------------------     S E C T I O N  3 -----------------------------------------
%******************************************************************************************
%##########################################################################################

\section{Description of the iterative approximate method}

Let us denote by \(P_{i}\) the hand side of \(i\)-th equation of the system \eqref{MainSystem}
\begin{equation}\label{Eqns}
 P_{i}=\sum_{j = 1}^n \int\limits_{\alpha_{j - 1}(t)}^{\alpha_j(t)} K_{ij}(t,s)G_{ij}(s,x_j(s))ds - f_i(t),
\end{equation}
and denote by \(X\) the set of desired functions, and by  \(F\)  the vector of source functions \(f_{i}\)
\[
     X =\left(
       \begin{array}{c}
         x_1(t) \\
         x_2(t) \\
         \vdots \\
         x_n(t)
       \end{array}
     \right), \quad
     F =\left(
       \begin{array}{c}
         f_1(t) \\
         f_2(t) \\
         \vdots \\
         f_n(t)
       \end{array}
     \right).
\]
We also introduce a matrix integral operator
\begin{equation}\label{Operator}
    P(X) \equiv \begin{pmatrix}
    P_{11} & P_{12} & \cdots & P_{1n} \\
    P_{21} & P_{22} & \cdots & P_{2n} \\
    \vdots & \vdots & \ddots & \vdots\\
    P_{n1} & P_{n2} & \cdots & P_{nn} \\
    \end{pmatrix}
    \begin{pmatrix}
      x_{1} \\
      x_{2} \\
      \vdots \\
      x_{n}
    \end{pmatrix} -
    \begin{pmatrix}
       f_{1} \\
       f_{2} \\
       \vdots \\
       f_{n}
   \end{pmatrix},
\end{equation}
which components are nonlinear integral operators
\[
 (P_{ij}x)(t) \equiv \int\limits_{\alpha_{j - 1}(t)}^{\alpha_j(t)} K_{ij}(t,s)G_{ij}(s,x(s))ds, \quad i,j=\overline{1,n}.
 \]
Let us rewrite system \eqref{MainSystem}  in the operator form
\begin{equation}\label{Eq21}
   P(X) = 0.
\end{equation}
To construct an iterative numerical method for solving  equation \eqref{Eq21},  the following modified Newton-Kantorovich scheme \cite{Kantorovich2} is employed
\begin{equation}\label{N-K}
   X^{m + 1} = X^m - [P'(X^0)]^{-1}(P(X^m)), m = 0,1,...,
\end{equation}
where \(  X^0 =
 \begin{pmatrix}
    x_{1}^{0}(t) \\
    x_{2}^{0}(t) \\
    \vdots \\
    x_{n}^{0}(t)
 \end{pmatrix}\) is
the is the set of initial approximation functions, and the solution of \eqref{Eq21} is defined as the limit
\(  X^{*} = \lim\limits_{m \to \infty} X^m  \).
Derivative \(P'(X^0)\) of operator \eqref{Operator} is as follows
\begin{equation}\label{Derivative}
   \begin{pmatrix}
     \frac{\partial P_{11}}{\partial x_1(t)} & \frac{\partial P_{12}}{\partial x_2(t)} & \cdots & \frac{\partial P_{1n}}{\partial x_n(t)} \\
     \frac{\partial P_{21}}{\partial x_1(t)} & \frac{\partial P_{22}}{\partial x_2(t)} & \cdots & \frac{\partial P_{2n}}{\partial x_n(t)} \\
     \vdots & \vdots & \ddots & \vdots\\
     \frac{\partial P_{n1}}{\partial x_1(t)} & \frac{\partial P_{n2}}{\partial x_2(t)} & \cdots & \frac{\partial P_{nn}}{\partial x_n(t)} \\
   \end{pmatrix}
   \begin{pmatrix}
      x_{1}^{0} \\
      x_{2}^{0} \\
      \vdots \\
      x_{n}^{0}
   \end{pmatrix}.
\end{equation}

Its components \(\left(\frac{\partial P_{ij}}{\partial x_j(t)}\right)(x_j^0),\; i,j=\overline{1,n},\) are defined as follows
\[
  \left(\frac{\partial P_{ij}}{\partial x_j(t)}\right)(x_j^0)(t) = \lim_{\omega \to 0} \frac{P_{ij}(x_j^0 + \omega x) - P_{ij}(x_j^0)}{\omega} =
\]
\[
 =\lim_{\omega \to 0} \frac{1}{\omega}  \int\limits_{\alpha_{j - 1}(t)}^{\alpha_j(t)} K_{ij}(t,s) [G_{ij}(s,x_j^0(s) + \omega x(s)) - G_{ij}(s,x_j^0(s))]ds.
\]
The limit transition under the sign of the integral implies:
\[
 \left(\frac{\partial P_{ij}}{\partial x_j(t)}\right)(x_j^0)(t) = \int\limits_{\alpha_{j-1}(t)}^{\alpha_j(t)}K_{ij}(t,s)G_{ij_x}(s,x_j^0(s))x(s)ds,
\]
where
\[
 G_{ij_x}(s,x_j^0(s)) = \frac{\partial G_{ij}(s,x_j(s))}{\partial x_j } \bigg|_{x = x_j^0}.
\]
On iterations with the number \(m\) we have an operator equation with respect to the correction
 \(\Delta X^{m + 1} = X^{m + 1} - X^m\), \(m=0,1,2,\ldots\)
\[
  P'(X^0) \Delta X^{m + 1} = -P(X^m),
\]
which in expanded form is equivalent to the following system of equations
\begin{equation*}
  \begin{cases}
    \sum\limits_{j = 1}^n \int\limits_{\alpha_{j - 1}(t)}^{\alpha_j(t)}K_{1j}(t,s) G_{1j}(s,x_j^0(s)) \Delta x_j^{m + 1}(s)ds = f_1(t) - \sum\limits_{j = 1}^n \int\limits_{\alpha_{j - 1}(t)}^{\alpha_j(t)}K_{1j}(t,s) G_{1j}(s,x_j^m(s))ds,\\
    \sum\limits_{j = 1}^n \int\limits_{\alpha_{j - 1}(t)}^{\alpha_j(t)}K_{2j}(t,s) G_{2j}(s,x_j^0(s)) \Delta x_j^{m + 1}(s)ds = f_2(t) - \sum\limits_{j = 1}^n \int\limits_{\alpha_{j - 1}(t)}^{\alpha_j(t)}K_{1j}(t,s) G_{2j}(s,x_j^m(s))ds,\\
    \vdots\\
    \sum\limits_{j = 1}^n \int\limits_{\alpha_{i - 1}(t)}^{\alpha_i(t)}K_{nj}(t,s) G_{nj}(s,x_j^0(s)) \Delta x_j^{m + 1}(s)ds = f_n(t) - \sum\limits_{j = 1}^n \int\limits_{\alpha_{j - 1}(t)}^{\alpha_j(t)}K_{1j}(t,s) G_{nj}(s,x_j^m(s))ds.\\
\end{cases}
\end{equation*}
The latter system, after a number of simplifications, is brought to mind
\begin{equation}\label{NKsystem}
  \begin{cases}
     \sum\limits_{j = 1}^n \int\limits_{\alpha_{j - 1}(t)}^{\alpha_j(t)}K_{1j}(t,s) G_{1j}(s,x_j^0(s)) x_j^{m + 1}(s)ds = \Psi^m_1(t),\\
     \sum\limits_{j = 1}^n \int\limits_{\alpha_{j - 1}(t)}^{\alpha_j(t)}K_{2j}(t,s) G_{2j}(s,x_j^0(s)) x_j^{m + 1}(s)ds = \Psi^m_2(t),\\
     \vdots\\
     \sum\limits_{j = 1}^n \int\limits_{\alpha_{i - 1}(t)}^{\alpha_i(t)}K_{nj}(t,s) G_{nj}(s,x_j^0(s)) x_j^{m + 1}(s)ds = \Psi^m_n(t).\\
  \end{cases}, \; m = 0, 1, 2,...,
\end{equation}
where
\[
  \Psi^m_i(t)=f_i(t) - \sum\limits_{j = 1}^n \int\limits_{\alpha_{j-1}(t)}^{\alpha_j(t)}K_{ij}(t,s)[G_{ij_x}(s,x^0_j(s))x^m_j(s) - G_{ij}(s,x^m_j(s))]ds.
\]

The system \eqref{NKsystem} is a system of linear Volterra integral equations with discontinuous kernels \(K_{ij}(t,s) G_{ij}(s,x^0_j(s))\), \(i=\overline{1,n}\), \(j=\overline{1,n}\) which do not change from iteration to iteration. A numerical method for solving similar systems of equations involved in the iterative process \eqref{N-K} is proposed below.

\subsection{The convergence theorem}

Let \(\overline{C}_{[0,T]}\) be a Banach vector-space of continuous vector-functions
\[
  \overline{C}_{[0,T]}=\bigl\{X=\left(x_1(t),x_2(t),\ldots,x_n(t)\right):\;x_i(t)\in C_{[0,T]}, \; i=\overline{1,n}\bigr\}
\]
equipped with the norm
\[
  \|X\|_{\overline{C}_{[0,T]}}=\max\limits_{i=\overline{1,n}}{\|x_i\|_{C_{[0,T]}}},  \text{ where }  \|x_i\|_{C_{[0,T]}}=\max\limits_{t\in [0,T]}|x_i(t)|,
  \; i=\overline{1,n}.
\]
The following theorem of convergence (based on the general theory proposed in the classical monograph \cite{Kantorovich2}) for iterative process \eqref{N-K} takes place:
%%------------------
\begin{theorem}
   Let the operator \(P\) has a continuous second derivative in the sphere
   \(\Omega_0\;(\|X-X^0\|_{\overline{C}_{[0,T]}}\leqslant r)\) and the following conditions
   hold:
   \begin{enumerate}
      \item Equation \eqref{NKsystem} has a unique solution for \(m=0\), i.e. there exists \(\Upsilon_0=[P'(X^0)]^{-1}\);
      \item \(\|\Delta X^1\|\leqslant\eta,\) where \(\Delta X^m=X^m-X^{m-1},\; m=0,1,\ldots.\);
      \item \(\|\Upsilon_0P''(X)\|\leqslant L,\;\;X\in\Omega_0\).
   \end{enumerate}
If also
\(
  h=L\eta<\frac12 \text{ and }
  \frac{1-\sqrt{1-2h}}{h}\eta\leqslant
  r\leqslant\frac{1+\sqrt{1-2h}}{h}\eta,
\)
then equation \eqref{Eq21} has a unique solution \(X^*\) in \(\Omega_0\), process \eqref{N-K} converges to \(X^*\), and the
velocity of convergence is estimated by the inequality
  \[\|X^*-X^m\|\leqslant \frac{\eta}{h}(1-\sqrt{1-2h})^{m+1},
  \;m=0,1,\ldots.\]
\end{theorem}

In order to prove this theorem one must show that equation \eqref{NKsystem} is uniquely solvable (including the case \(m=0\)),
i.e. condition 1 of the theorem holds. The solvability of linear systems of such Volterra integral equations was discussed in \cite{Sidorov-2013}.  The further argumentation is performed under conditions provided the uniqueness of solution for such systems (see \cite{Sidorov-2013} for details).

Then the boundedness of the second derivative \(\bigl[P''(X^0)\bigr](X)\) should be verified
for estimating the constant \(L\) in condition 3. One can verify that the necessary condition for the second derivative \(\bigl[P''(X^0)\bigr](X) \) to be bounded is a
differentiability of the initial approximation \(X^0\) as well as the functions \(K_{ij}\) with respect to second variable and the functions \(G_{ij}\) with respect to both variables.

%##########################################################################################
%******************************************************************************************
%---------------------------     S E C T I O N  3 -----------------------------------------
%******************************************************************************************
%##########################################################################################

\section{Discretization of the system of linear integral equations}
\subsection{Problem formulation}

In this paragraph, we propose a numerical method for solving systems of linear integral equations arising at each step of the iterative process \eqref{N-K}.

Let us consider the following system of linear integral equations of the first kind
\begin{equation}\label{LinSystem}
   \begin{cases}
      \sum\limits_{j = 1}^n \int\limits_{\alpha_{j-1}(t)}^{\alpha_j(t)} h_{1}(t,s)x_j(s)ds = f_1(t),\\
      \sum\limits_{j = 1}^n \int\limits_{\alpha_{j-1}(t)}^{\alpha_j(t)} h_{2}(t,s)x_j(s)ds = f_2(t),\\
      \vdots\\
      \sum\limits_{j = 1}^n \int\limits_{\alpha_{j-1}(t)}^{\alpha_j(t)} h_{n}(t,s)x_j(s)ds = f_n(t).\\
   \end{cases},
\end{equation}
 functions \(h_{i}(t, s), i=\overline{1,n},\)  undergo a finite discontinuity on the curves \(\alpha_j(t), j =\overline{1,n-1}\), and
\begin{equation*}
   h_{i}(t, s) =
  \begin{cases}
     K_{i1}:=K_{i1}(t,s)G_{i1}(s,x_j^0(s)),&\text{ if } \alpha_0(t) < s \le \alpha_1(t);  \\
     K_{i2}:=K_{i2}(t,s)G_{i2}(s,x_j^0(s)), &\text{ if } \alpha_1(t) < s \le \alpha_2(t); \\
     \vdots\\
     K_{in}:=K_{in}(t,s)G_{in}(s,x_j^0(s)), &\text{ if } \alpha_{n-1}(t) < s \le \alpha_n(t).
\end{cases}
\end{equation*}
Here  \(\alpha_0(t) \equiv 0,\; \alpha_0(t) <  \alpha_1(t) < ... < \alpha_n(t) \equiv t, \; f(0) = 0\).
The new kernels \(K_{ij}(t,s)\) and right-hand side \(f_i(t)\) are continuous and sufficiently smooth functions introduced to simplify the exposition. The functions \(\alpha_j(t) \in C^1 [0,T]\) are non-decreasing.
Also \(\alpha_1'(0) \le \alpha_2'(0) \le ... \le \alpha_{n-1}'(0) < 1\). A detailed theoretical study of equations of this type is carried out in the book \cite{Sidorov-book}.

%********************************************************************
%********************************************************************
%********************************************************************

\subsection{Piecewise constant approximation}
\hspace * {\parindent}
In order to construct a numerical solution of \eqref{LinSystem} on the segment \([0,T]\) (under the conditions of existence of a single continuous solution), we introduce a grid of nodes (not necessarily uniform)
\[
   0 = t_0 < t_1 < t_2 < ... < t_N = t,\; h = \max\limits_{j = \overline{1, N}}(t_j - t_{j-1}) = O\left(N^{-1}\right).
\]
Approximate solution of the system \eqref{LinSystem} will be sought in the form of a vector of piecewise constant functions
\begin{equation*}
   X^N (t) =
   \begin{pmatrix}
     \sum\limits_{j = 1}^N x_1^j \delta_j(t) \\
     \sum\limits_{j = 1}^N x_2^j \delta_j(t)\\
     \vdots \\
     \sum\limits_{j = 1}^N x_n^j \delta_j(t)
   \end{pmatrix},\; t \in (0,T], \;
   \delta_j (t) =
   \begin{cases}
     1, &\text{when } t \in \Delta_j = (t_{j-1}, t_j], \\
     0, &\text{when } t \notin \Delta_j.
   \end{cases},
\end{equation*}
with uncertain coefficients so far \(x_i^j, \; i = \overline{1,n}, \; j = \overline{1,N}.\)

To determine the values \(x_i^0 = x_i(0), \; i = \overline{1,n},\) we differentiate both parts of each equation of \eqref{LinSystem} by \(t\)
\[
f _i(t) = \sum\limits_{j = 1}^n \Bigg( \int\limits_{\alpha_{j-1}(t)}^{\alpha_j(t)} \frac{\partial K_{ij}(t,s)}{\partial t} x_j(s) ds +\alpha_j'(t)K_{ij}(t,\alpha_j(t))x_j(\alpha_j(t)) - \alpha_{j-1}'(t)K_{ij}(t,\alpha_{j-1}(t))x_j(\alpha_{i-1}(t))\Bigg).
\]

Using the latter ratio, we come to the necessity of solving a system of equations
\begin{equation}\label{SLAU}
  \begin{cases}
     \sum\limits_{j=1}^n K_{1j}(0,0) [\alpha_j'(0) - \alpha_{j-1}'(0)] x_j^0 = f ' _1(0),\\
     \sum\limits_{j=1}^n K_{2j}(0,0) [\alpha_j'(0) - \alpha_{j-1}'(0)] x_j^0 = f ' _2(0),\\
     \vdots\\
     \sum\limits_{j=1}^n K_{nj}(0,0) [\alpha_j'(0) - \alpha_{j-1}'(0)] x_j^0 = f'_n(0).
  \end{cases}
\end{equation}
with respect to values \(x_j^0, \; j = \overline{1,n}\). Here it is assumed that the components of the equation \eqref{LinSystem} are such that the system of linear algebraic equations \eqref{SLAU} has a unique nontrivial solution.

Next, we introduce the designation \(f_i^k = f_i (t_k),\; I=1,\ldots,n, \; k = 1,\ldots,N\). To determine the coefficients \(x_j^1\), write down each equation of the system \eqref{LinSystem} at \(t = t_1\).
\begin{equation}\label{Eq41}
   \sum\limits_{j=1}^n \int\limits_{\alpha_{j-1}(t_1)}^{\alpha_j(t_1)} K_{ij}(t_1,s)x_j(s)ds = f_i^1.
\end{equation}
Since at this step the lengths of all integration segments \(\alpha_j(t_1) - \alpha_{j-1}(t_1)\) in \eqref{Eq41} do not exceed the maximum grid step \(h\), and the components of the approximate solution take the values \(x_j^1\), \(j=\overline{1,n}\), then applying the quadrature formula of the middle rectangles, we have a system of
\begin{equation*}
\begin{cases}
\sum\limits_{j=1}^n (\alpha_j (t_1) - \alpha_{j-1} (t_1)) K_{1j}(t_1,\frac{\alpha_j(t_1) + \alpha_{j-1}(t_1)} {2}) x_j^1 = f_1^1,\\
\sum\limits_{j=1}^n (\alpha_j (t_1) - \alpha_{j-1} (t_1)) K_{2j}(t_1,\frac{\alpha_j (t_1) + \alpha_{j-1} (t_1)}{2}) x_j^1 = f_2^1,\\
\vdots\\
\sum\limits_{j=1}^n (\alpha_j(t_1) - \alpha_{j-1}(t_1))K_{nj}(t_1,\frac{\alpha_j(t_1) + \alpha_{j-1}(t_1)}{2}) x_j^1 = f_n^1.\\
\end{cases},
\end{equation*}
of which shall be determined by the value of \(x_j^1,\;j=\overline{1,n}\).

Let further \(l_j^K\) denote the number of the grid's segment, which gets the value \(\alpha_{j} (t_k)\), that is, the condition \(t_{l_j^k-1}\leq \alpha_{j}(t_k) \leq t_{l_j^k}\) is satisfied.

Assume now that the values \(x_j^2, x_j^3, ..., x_j^{l_j^K-1},\;K=\overline{1,n}\) are already found. Let us transform system \eqref{LinSystem} and require the following equalities satisfaction at the point \(t = t_k\)
\begin{equation*}
   \begin{cases}
      \sum\limits_{j=1}^n\int\limits_{t_{l_j^k-1}}^{\alpha_{j}(t_{k})} K_{1j}(t_{k},s)x_j(s)ds = f_1^k - \sum\limits_{j=1}^n\int\limits_{\alpha_{j-1}(t_k)}^{t_{l_j^k-1}} K_{1j}(t_{k},s)x_j^N(s)ds,\\
      \sum\limits_{j=1}^n\int\limits_{t_{l_j^k-1}}^{\alpha_{j}(t_{k})} K_{2j}(t_{k},s)x_j(s)ds = f_2^k - \sum\limits_{j=1}^n\int\limits_{\alpha_{j-1}(t_k)}^{t_{l_j^k-1}} K_{2j}(t_{k},s)x_j^N(s)ds,\\
      \vdots\\
      \sum\limits_{j=1}^n\int\limits_{t_{l_j^k-1}}^{\alpha_{j}(t_{k})} K_{nj}(t_{k},s)x_j(s)ds = f_n^k - \sum\limits_{j=1}^n\int\limits_{\alpha_{j-1}(t_k)}^{t_{l_j^k-1}} K_{nj}(t_{k},s)x_j^N(s)ds.\\
\end{cases}
\end{equation*}
Given approximation leads to the following  system
\begin{equation}\label{Eq42}
  \begin{cases}
    \sum\limits_{j=1}^n x_j^{l_j^k}\int\limits_{t_{l_j^k-1}}^{\alpha_{j}(t_{k})} K_{1j}(t_{k},s)ds = f_1^k - \sum\limits_{j=1}^n\int\limits_{\alpha_{j-1}(t_k)}^{t_{l_j^k-1}} K_{1j}(t_{k},s)x_j^N(s)ds,\\
    \sum\limits_{j=1}^n x_j^{l_j^k}\int\limits_{t_{l_j^k-1}}^{\alpha_{j}(t_{k})} K_{2j}(t_{k},s)ds = f_2^k - \sum\limits_{j=1}^n\int\limits_{\alpha_{j-1}(t_k)}^{t_{l_j^k-1}} K_{2j}(t_{k},s)x_j^N(s)ds,\\
    \vdots\\
    \sum\limits_{j=1}^n x_j^{l_j^k}\int\limits_{t_{l_j^k-1}}^{\alpha_{j}(t_{k})} K_{nj}(t_{k},s)ds = f_n^k - \sum\limits_{j=1}^n\int\limits_{\alpha_{j-1}(t_k)}^{t_{l_j^k-1}} K_{nj}(t_{k},s)x_j^N(s)ds.\\
  \end{cases}
\end{equation}

 For the integrals calculation in \eqref{Eq42} compound midpoint rectangle quadrature constructed on the auxiliary grid of nodes bound at each concrete value \(N\) to discontinuity curve \(\alpha_i(t)\) of kernel \(K(t, s)\) are used.

It is not difficult to see such an approximation method has the first order of accuracy and the following error estimate takes place
\[
  \varepsilon_N = \left\|X - X_N \right\|_{C[0,T]} = {\mathcal O}\left(\frac{1}{N}\right).
\]

\textbf{Remark.} \textit{It should be noted that the proposed numerical approach also allows a more accurate approximation of the solution. In particular, when using piecewise linear approximation, the order of accuracy increases by one. However such discretized systems may have nonunique solutions (due to the discontinuities \(\alpha_i(t), \; i=\overline{1,n},\)  they can be underdetermined or have many solutions). Taking into account this fact we suggest new alternative numerical method below.}

%******************************************************************************************
%******************************************************************************************
%******************************************************************************************

\subsection{Polynomial collocation method}\label{Section-LinearCase}
\subsubsection{Description of the problem}
In this section we deal with numerical solution of the following system of VIE
\begin{equation}\label{MainEq}
    \int\limits_{0}^{t} H(t,s)X(s)ds = F(t),
\end{equation}
where
\[
  H(t,s)=
  \left(\begin{array}{c}
    h_1(t,s) \\
    h_2(t,s) \\
    \vdots \\
    h_n(t,s)
  \end{array}\right), \quad
  F(t)=
  \left(\begin{array}{c}
    f_1(t) \\
    f_2(t) \\
    \vdots \\
    f_n(t)
  \end{array}\right).
\]
Here kernels \(h_{i}(t, s), i=\overline{1,n},\) have just discontinuities on the lines \(\alpha_j(t), j =\overline{1,n-1}\) and kernels defined as follows
\begin{equation*}
   h_{i}(t, s) =
  \begin{cases}
     K_{i1}(t,s), &\text{ if } \alpha_0(t) < s \le \alpha_1(t);  \\
     K_{i2}(t,s), &\text{ if } \alpha_1(t) < s \le \alpha_2(t); \\
     \vdots\\
     K_{in}(t,s), &\text{ if } \alpha_{n-1}(t) < s \le \alpha_n(t).
\end{cases}
\end{equation*}
Components \(x_i(t),\; i=\overline{1,n},\) of desired function \(X(t)\)
are defined in the domains \((\alpha_{i-1}(t), \alpha_i(t)]\), i.e.
\[
  X(t)=x_i(t) \text{ for } t\in \alpha_{i-1}(t) < s \le \alpha_i(t).
\]

Here \(\alpha_0(t) \equiv 0,\; \alpha_0(t) <  \alpha_1(t) < ... < \alpha_n(t) \equiv t, \) for \(t\in(0,T)\), \(\alpha_i(0)=f_i(0) = 0, \; i=\overline{1,n}\).
It is assumed that functions \(K_{ij}(t,s)\) and \(f_i(t), \;i,j=\overline{1,n},\) are
continuous and smooth sufficiently. Functions \(\alpha_j(t) \in  C^1 [0,T]\) and nondecreasing, \(\alpha_1'(0) \le \alpha_2'(0) \le ... \le \alpha_{n-1}'(0) < 1\).

Let us rewrite the formal notation \eqref{MainEq} as follows

\begin{equation}\label{LinSystem1}
   \begin{cases}
      \sum\limits_{j = 1}^n \int\limits_{\alpha_{j-1}(t)}^{\alpha_j(t)} K_{1j}(t,s)x_j(s)ds = f_1(t),\\
      \sum\limits_{j = 1}^n \int\limits_{\alpha_{j-1}(t)}^{\alpha_j(t)} K_{2j}(t,s)x_j(s)ds = f_2(t),\\
      \vdots\\
      \sum\limits_{j = 1}^n \int\limits_{\alpha_{j-1}(t)}^{\alpha_j(t)} K_{nj}(t,s)x_j(s)ds = f_n(t).\\
   \end{cases}.
\end{equation}

%---------------------------------------------------------------------
%----------------------  S U B S E C T I O N S -----------------------
%---------------------------------------------------------------------
\subsubsection{Collocation}
\hspace*{\parindent}
For construction of numerical solution of system \eqref{LinSystem1} on \([0,T]\)
(in conditions of unique continuous solution existence) let  us introduce the mesh (not necessarily uniform in general)
\[
   0 = t_0 < t_1 < t_2 < ... < t_m = T,\;\,\, h = \max\limits_{j = \overline{1,m}}(t_j - t_{j-1}) = {\mathcal O} \left(m^{-1}\right).
\]
An approximate solution of the system \eqref{LinSystem1} we search as
\begin{equation}\label{Xm}
   X^m(t) =
   \begin{pmatrix}
     \sum\limits_{j = 0}^m A_{1j} t^j \\
     \sum\limits_{j = 0}^m A_{2j} t^j\\
     \vdots \\
     \sum\limits_{j = 0}^m A_{nj} t^j
   \end{pmatrix}, \; t \in (0,T],
\end{equation}
with unknown coefficients \(A_{ij}, \; i = \overline{1,n}, \; j = \overline{0,m}.\)

To define the coefficients \(A_{i0} = x_i(0), \; i = \overline{1,n},\) we
differetiate  eq. \eqref{LinSystem1} with respect to \(t\)
\[
f'_i(t) = \sum\limits_{j = 1}^n \Bigg( \int\limits_{\alpha_{j-1}(t)}^{\alpha_j(t)} \frac{\partial K_{ij}(t,s)}{\partial t} x_j(s) ds +\alpha_j'(t)K_{ij}(t,\alpha_j(t))x_j(\alpha_j(t)) - \alpha_{j-1}'(t)K_{ij}(t,\alpha_{j-1}(t))x_j(\alpha_{i-1}(t))\Bigg).
\]
Using the latter equation we obtain the following system
\begin{equation}\label{SLAU}
  \begin{cases}
     \sum\limits_{j=1}^n K_{1j}(0,0) [\alpha_j'(0) - \alpha_{j-1}'(0)] A_{j0} = f'_1(0),\\
     \sum\limits_{j=1}^n K_{2j}(0,0) [\alpha_j'(0) - \alpha_{j-1}'(0)] A_{j0} = f'_2(0),\\
     \vdots\\
     \sum\limits_{j=1}^n K_{nj}(0,0) [\alpha_j'(0) - \alpha_{j-1}'(0)] A_{j0} = f'_n(0).
  \end{cases}
\end{equation}
with respect to  \(A_{j0}, \; j = \overline{1,n}\). Here we assume that components of
eq. \eqref{LinSystem1} such as SLAE \eqref{SLAU}
enjoys unique nontrivial solution.

To define the coefficients \(A_{ij}\) each eq. of the system \eqref{LinSystem1}
we write in the point \(t = t_k\).
\begin{equation}\label{Eq41}
   \sum\limits_{j=1}^n \int\limits_{\alpha_{j-1}(t_k)}^{\alpha_j(t_k)} K_{ij}(t_k,s)x_j(s)ds = f_{i}(t_k).
\end{equation}
Let us require the conversion of eq. \eqref{Eq41} (considering \eqref{Xm})  into the following identity:
\begin{equation}\label{Eq42}
   \sum\limits_{j=1}^n \int\limits_{\alpha_{j-1}(t_k)}^{\alpha_j(t_k)} K_{ij}(t_k,s)\left(\sum\limits_{l=0}^m A_{jl} s^l\right) ds = f_{i}(t_k).
\end{equation}
Because coefficients \(A_{j0}, \; j = \overline{1,n} \) already known, we can write
\begin{equation}\label{Eq42}
   \sum\limits_{j=1}^n \sum\limits_{l=1}^m  A_{jl} \int\limits_{\alpha_{j-1}(t_k)}^{\alpha_j(t_k)} K_{ij}(t_k,s) s^l ds = f_{i}(t_k)-\sum\limits_{j=1}^n A_{j0}\int\limits_{\alpha_{j-1}(t_k)}^{\alpha_j(t_k)} K_{ij}(t_k,s) ds.
\end{equation}
Let us introduce the following notations
\begin{equation}\label{F}
F_{ik} := f_{i}(t_k)-\sum\limits_{j=1}^n A_{j0}\int\limits_{\alpha_{j-1}(t_k)}^{\alpha_j(t_k)} K_{ij}(t_k,s) ds; \; i=\overline{1,n}, k = \overline{1,m},
\end{equation}
\begin{equation}\label{C}
C_{ikjl} := \int\limits_{\alpha_{j-1}(t_k)}^{\alpha_j(t_k)} K_{ij}(t_k,s) s^l ds; \; i=\overline{1,n}, k = \overline{1,m}, j = \overline{1,n}, l = \overline{1,m};
\end{equation}
\begin{equation}\label{Cikjl}
{\tiny
C:=
\begin{pmatrix}
    \begin{tabular}{cccc}
        \(C_{1111}\) & \(C_{1112}\)  &\ldots  & \(C_{111m}\)\\
        \(C_{1211}\) &\(C_{1212}\)   &\ldots  & \(C_{121m}\)\\
        \ldots       &\ldots         & \ldots & \ldots\\
        \(C_{1m11}\) &\(C_{1m12}\)   &\ldots  & \(C_{1m1m}\)
    \end{tabular}
     &
    \begin{tabular}{cccc}
        \(C_{1121}\) & \(C_{1122}\) & \ldots  & \(C_{112m}\)\\
        \(C_{1121}\) &\(C_{1222}\)  &\ldots   &\(C_{122m}\)\\
        \ldots       &\ldots        &\ldots   &\ldots\\
        \(C_{1m21}\) &\(C_{1m22}\)  &\ldots   &\(C_{1m2m}\)
    \end{tabular}

    &

    \begin{tabular}{cccc}
      \ldots\\
      \ldots\\
      \ldots\\
      \ldots\\
    \end{tabular}

    &

    \begin{tabular}{cccc}
        \(C_{11n1}\) & \(C_{11n2}\) & \ldots  & \(C_{11nm}\)\\
        \(C_{12n1}\) &\(C_{12n2}\)  &\ldots   &\(C_{12nm}\)\\
        \ldots       &\ldots        &\ldots   &\ldots\\
        \(C_{1mn1}\) &\(C_{1mn2}\)  &\ldots   &\(C_{1mnm}\)
    \end{tabular}\\
 %------------------
        \begin{tabular}{cccc}
        \(C_{2111}\) & \(C_{2112}\)  &\ldots  & \(C_{211m}\)\\
        \(C_{2211}\) &\(C_{2212}\)   &\ldots  & \(C_{221m}\)\\
        \ldots       &\ldots         & \ldots & \ldots\\
        \(C_{2m11}\) &\(C_{2m12}\)   &\ldots  & \(C_{2m1m}\)
    \end{tabular}
     &
    \begin{tabular}{cccc}
        \(C_{2121}\) & \(C_{2122}\) & \ldots  & \(C_{212m}\)\\
        \(C_{2121}\) &\(C_{2222}\)  &\ldots   &\(C_{222m}\)\\
        \ldots       &\ldots        &\ldots   &\ldots\\
        \(C_{2m21}\) &\(C_{2m22}\)  &\ldots   &\(C_{2m2m}\)
    \end{tabular}
    &
    \begin{tabular}{cccc}
      \ldots\\
      \ldots\\
      \ldots\\
      \ldots\\
    \end{tabular}
    &
    \begin{tabular}{cccc}
        \(C_{21n1}\) & \(C_{21n2}\) & \ldots  & \(C_{21nm}\)\\
        \(C_{22n1}\) &\(C_{22n2}\)  &\ldots   &\(C_{22nm}\)\\
        \ldots       &\ldots        &\ldots   &\ldots\\
        \(C_{2mn1}\) &\(C_{2mn2}\)  &\ldots   &\(C_{2mnm}\)
    \end{tabular}\\
  %------------------
    \begin{tabular}{cccc}
        \ldots       &\ldots        &\ldots   &\ldots\\
        \ldots       &\ldots        &\ldots   &\ldots\\
        \ldots       &\ldots        &\ldots   &\ldots\\
        \ldots       &\ldots        &\ldots   &\ldots\\
    \end{tabular}
    &
    \begin{tabular}{cccc}
        \ldots       &\ldots        &\ldots   &\ldots\\
        \ldots       &\ldots        &\ldots   &\ldots\\
        \ldots       &\ldots        &\ldots   &\ldots\\
        \ldots       &\ldots        &\ldots   &\ldots\\
    \end{tabular}
    &
    \begin{tabular}{cccc}
      \ldots\\
      \ldots\\
      \ldots\\
      \ldots\\
    \end{tabular}
    &
    \begin{tabular}{cccc}
        \ldots       &\ldots        &\ldots   &\ldots\\
        \ldots       &\ldots        &\ldots   &\ldots\\
        \ldots       &\ldots        &\ldots   &\ldots\\
        \ldots       &\ldots        &\ldots   &\ldots\\
    \end{tabular}\\
 %------------------
        \begin{tabular}{cccc}
        \(C_{n111}\) & \(C_{n112}\)  &\ldots  & \(C_{n11m}\)\\
        \(C_{n211}\) &\(C_{n212}\)   &\ldots  & \(C_{n21m}\)\\
        \ldots       &\ldots         & \ldots & \ldots\\
        \(C_{nm11}\) &\(C_{nm12}\)   &\ldots  & \(C_{nm1m}\)
    \end{tabular}
     &
    \begin{tabular}{cccc}
        \(C_{n121}\) & \(C_{n122}\) & \ldots  & \(C_{n12m}\)\\
        \(C_{n121}\) &\(C_{n222}\)  &\ldots   &\(C_{n22m}\)\\
        \ldots       &\ldots        &\ldots   &\ldots\\
        \(C_{nm21}\) &\(C_{nm22}\)  &\ldots   &\(C_{nm2m}\)
    \end{tabular}
    &
    \begin{tabular}{cccc}
       \ldots\\
       \ldots\\
       \ldots\\
       \ldots\\
    \end{tabular}
    &
    \begin{tabular}{cccc}
        \(C_{n1n1}\) & \(C_{n1n2}\) & \ldots  & \(C_{n1nm}\)\\
        \(C_{n2n1}\) &\(C_{n2n2}\)  &\ldots   &\(C_{n2nm}\)\\
        \ldots       &\ldots        &\ldots   &\ldots\\
        \(C_{nmn1}\) &\(C_{nmn2}\)  &\ldots   &\(C_{nmnm}\)
    \end{tabular}\\

\end{pmatrix},
}
\end{equation}

\begin{equation}\label{AF}
A:=
\begin{pmatrix}
 A_{11}\\
 A_{12}\\
\ldots\\
 A_{1m}\\
 A_{21}\\
 A_{22}\\
 \ldots\\
 A_{2m}\\
\ldots \\
 A_{n1}\\
 A_{n2}\\
 \ldots\\
 A_{nm}\\
\end{pmatrix}, F:=
\begin{pmatrix}
 F_{11}\\
 F_{12}\\
\ldots\\
 F_{1m}\\
 F_{21}\\
 F_{22}\\
 \ldots\\
 F_{2m}\\
\ldots \\
 F_{n1}\\
 F_{n2}\\
 \ldots\\
 F_{nm}\\
\end{pmatrix}.
\end{equation}
As result in order to define the coeffients \(A_{ij}\)
one must solve the folowing SLAE of  dimension \((n\cdot m)\times(n\cdot m)\)
\begin{equation}\label{MatrixSystem}
C \cdot A=F.
\end{equation}
Elements of the matices \eqref{Cikjl} and \eqref{AF} can be represented as follows
\begin{equation}
C_{[(i-1)\cdot m + k][(j-1)\cdot m + l]}=C_{ikjl};
\end{equation}
\begin{equation}
A_{[(j-1)\cdot m + l]}=A_{jl};
\end{equation}
\begin{equation}
F_{[(i-1)\cdot m + k]}=F_{ik};
\end{equation}
where \(i=\overline{1,n}, \; k = \overline{1,m}, \; j = \overline{1,n}, \; l = \overline{1,m}\).

Here for integrals calculation in \eqref{F}, \eqref{C} one must
employ compound midpoint quadrature rule constructed on the auxiliary mesh
linked to the lines  \(\alpha_i(t)\) of the kernels \(h_i(t,s), \; i=\overline{1,n}\) discontinueties.

Finally, as solution to the system \eqref{MatrixSystem} the coeffients \(A_{ij}\) will be determined and the polynomial approximation of the components of the approximate solution \(X^m(t)\)
is derived.

%##########################################################################################
%******************************************************************************************
%---------------------------    NUMERICAL RESULTS -----------------------------------------
%******************************************************************************************
%##########################################################################################

\section{Numerical results}

\subsection{System of linear integral equations}

In order to illustrate the convergence of proposed in section 3.3 algorithm we consider two system of equations.\\

\textbf{Example 1}  (\(n=2\)).
 \begin{equation}\label{Model-01}
 \begin{cases}
    \int\limits_{0}^{\frac{t}{2}} (1+t+s) x_1(s)ds + \int\limits_{\frac{t}{2}}^{t} x_2(s)ds= f_1(t)\\
    \int\limits_{0}^{\frac{t}{2}}  (1+t-s)x_1(s)ds - \int\limits_{\frac{t}{2}}^{t} x_2(s)ds= f_2(t)
 \end{cases}, \quad t \in [0,2],
 \end{equation}
 where
\[
 f_1(t) = 1+\sin\left(\frac{t}{2}\right) + \frac{3t}{2}\cdot \sin\left(\frac{t}{2}\right) + 2\cos\left(\frac{t}{2}\right) - \cos(t),
\]
\[
  f_2(t) = 1 + \sin\left(\frac{t}{2}\right) + \frac{t}{2} \sin\left(\frac{t}{2}\right) - 2\cos\left(\frac{t}{2}\right) + \cos(t).
\]

Exact solution of the system \eqref{Model-01} is \( x_1(t) = \cos(t),\; x_2(t) = \sin(t)\).
Numerical solution of the system \eqref{Model-01} is given in Tab. 1, where
\(m\) is order of approximating polynomial, \(\varepsilon_i=\max\limits_{t \in [0,a_i(T)]}| x^{*}_i(t) - x_i(t)|\),  \(t_{max}^i\)  is the greatest error point  \(\varepsilon_i,\; i=1,\ldots,n\), \(\varepsilon = \sqrt{\sum\limits_{i=1}^{n}\varepsilon_i^2} \).
%---------------------
\begin{table}[h]
 \centering
\begin{tabular}{|c|c|c|c|c|c|}
\hline
    \(m\) & \(\varepsilon_1\)      &\(t_{max}^1\) & \(\varepsilon_2\)        & \(t_{max}^2\)  & \(\varepsilon\) \\
\hline
    2  &\(9.82294\cdot 10^{-3}\)& 1.0          &  \(6.72940\cdot 10^{-2}\)& 2.0            & \(6.80072\cdot 10^{-2}\)\\
\hline
    3  &\(1.60472\cdot 10^{-3}\)&1.0           &  \(2.35676\cdot 10^{-2}\)& 2.0            & \(2.36222\cdot 10^{-2}\)\\
\hline
    5  &\(6.67315\cdot 10^{-6}\)&1.0           &  \(3.95344\cdot 10^{-4}\)& 2.0            & \(3.95400\cdot 10^{-4}\)\\
\hline
    8  &\(1.72968\cdot 10^{-8}\)&1.0           &  \(1.80165\cdot 10^{-7}\)& 2.0            & \(1.80994\cdot 10^{-7}\)\\
\hline
    12 &\(1.31014\cdot 10^{-8}\)&1.0           &  \(4.12674\cdot 10^{-9}\)& 2.0            & \(1.37360\cdot 10^{-8}\)\\
\hline
    15 &\(2.79384\cdot 10^{-8}\)&1.0           &  \(1.75046\cdot 10^{-8}\)& 2.0            & \(3.29692\cdot 10^{-8}\)\\
\hline
\end{tabular}
  \caption{Results for the problem \eqref{Model-01}}
\end{table}
%---------------------

\textbf{Example 2 } (\(n=3\)).

Let us consider the following system of equations
\begin{equation}\label{Model-02}
 \begin{cases}
\int\limits_{0}^{\frac{t}{3}} (1+t+s) x_1(s)ds + \int\limits_{\frac{t}{3}}^{\frac{2t}{3}}  x_2(s)ds+ \int\limits_{\frac{2t}{3}}^{t} (1+s) x_3(s)ds= f_1(t)\\
\int\limits_{0}^{\frac{t}{3}} (1+t-s) x_1(s)ds - \int\limits_{\frac{t}{3}}^{\frac{2t}{3}}  x_2(s)ds+ \int\limits_{\frac{2t}{3}}^{t} (1-t) x_3(s)ds= f_2(t)\\
\int\limits_{0}^{\frac{t}{3}} (1+\frac{t}{5}+s) x_1(s)ds + \int\limits_{\frac{t}{3}}^{\frac{2t}{3}}( 1+t-s) x_2(s)ds- \int\limits_{\frac{2t}{3}}^{t} (1+s) x_3(s)ds= f_3(t)
\end{cases}, \quad t \in [0,2].
\end{equation}
where \(f_i(t), \; i=1,2,3,\) are selected such as exact solution is as follows
\[   x_1(t) = \cos(t),\; x_2(t) = \sin(t),\; x_3 =\frac{1}{4}\sin(t). \]

In tab. 2 the numerical solution of system  \eqref{Model-02} is presented.

\begin{table}[h]
 \centering
\begin{tabular}{|c|c|c|c|c|c|c|c|}
\hline
     \(m\)  & \(\varepsilon_1\)  &\(t_{max}^1\)  & \(\varepsilon_2\)  & \(t_{max}^2\)  & \(\varepsilon_3\) & \(t_{max}^3\)  & \(\varepsilon\)\\
\hline
    2  &\(1.75798\cdot 10^{-3}\)   & 0.6(6)   &  \(3.30114\cdot 10^{-3}\) & 1.3(3)  &  \(6.72940\cdot 10^{-2}\) & 2.0  & \(3.44752\cdot 10^{-2}\)\\
\hline
    3  &\(2.75534\cdot 10^{-4}\)   & 0.6(6)   &  \(3.77383\cdot 10^{-3}\) & 1.3(3)  &  \(4.51700\cdot 10^{-3}\) & 2.0  & \(5.89245\cdot 10^{-3}\)\\
\hline
    5  &\(1.68420\cdot 10^{-6}\)   & 0.6(6)   &  \(2.89243\cdot 10^{-5}\) & 1.3(3)  &  \(9.1497\cdot 10^{-5}\)  & 2.0  & \(9.59747\cdot 10^{-5}\)\\
\hline
   8   &\(7.86707\cdot 10^{-10}\)  & 0.6(6)   &  \(6.86599\cdot 10^{-9}\) & 1.3(3)  &  \(4.15579\cdot 10^{-8}\) & 2.0  & \(4.21286\cdot 10^{-8}\)\\
\hline
  12   &\(1.33907\cdot 10^{-9}\)   & 0.6(6)   &  \(1.55675\cdot 10^{-9}\) & 1.3(3)  &  \(1.97100\cdot 10^{-10}\) & 2.0 & \(2.06287\cdot 10^{-9}\)\\
\hline
  15   &\(1.69313\cdot 10^{-9}\)   & 0.6(6)   &  \(1.00072\cdot 10^{-9}\) & 1.3(3)  &  \(1.43989\cdot 10^{-9}\)  & 2.0 & \(2.43750\cdot 10^{-9}\)\\
\hline
\end{tabular}
 \caption{Results for problem \eqref{Model-02}}
\end{table}

Therefore, we suggested the efficient method for solution of  systems \eqref{MainEq} for relatively small \(T\). The polynomial spline-approximation can be used to find the solution
larger interval.

%$$$$$$$$$$$$$$$$$$$$$$$$$$$$$$$$$$$$$$$$$$$$$$$$$$$$$$$$$$$$$$$$$$$$$$$$$$$$$$$$$$$$$$$$$$

\subsection{Nonlinear equations}

Let us first consider the following scalar VIE:
\begin{equation}\label{NonlinearScalar}
\int\limits_{0}^{t}K(t,s,x(s))ds=\frac{1}{3}t^3+\frac{123}{192}t^4+\frac{1}{160}t^5+\frac{17}{1920}t^6,
\end{equation}
where
\[
K(t,s,x(s))=
 \begin{cases}
   (1+t+s)(x(s)+x^2(s)), &\text{$0\le s \le \frac{1}{2}t$};\\
   (1+2t)x(s), &\text{$\frac{1}{2}t<s\le t$}.
 \end{cases}
\]
Here \(x^*(t)=t^2\),  \(t \in [0,1]\) is the exact solution to equation \eqref{NonlinearScalar}.

Results concerning solution of this equation using Newton-Kantorovich method (combined with direct discretization of linear equations described in Section 3.2) are shown in Tab.~\ref{Num-1}, where the following notations were used:
\(h\) is the mesh step,  \(m\) is iterations number of Newton-Kantorovich method, \(\varepsilon=\max\limits_{t\in[0,1]}|x_N(t)-x^*(t)|\).

\begin{table}[h!]\label{Num-1}
\centering
\begin{tabular}{|*{6}{c|}}\hline
\multicolumn{6}{|c|}{Piece-wise constant approximation} \\ \hline
h  & $\frac{1}{32}$ & $\frac{1}{64}$ & $\frac{1}{128}$ & $\frac{1}{256}$ & $\frac{1}{512}$  \\ \hline
m  & 5 & 5 & 5 & 5 & 5   \\ \hline
$\varepsilon$ & 0.0286877  & 0.0152708  & 0.00730057       & 0.0044031      & 0.00386043       \\ \hline
\multicolumn{6}{|c|}{Piece-wise linear approximation} \\ \hline
h  & $\frac{1}{32}$ & $\frac{1}{64}$ & $\frac{1}{128}$ & $\frac{1}{256}$ & $\frac{1}{512}$   \\ \hline
m  & 6 & 9 & 9 & 9 & 10   \\ \hline
$\varepsilon$ & 0.001496942  & 0.0008320731  & 0.0004604894       & 0.0002356853      & 0.0001186141        \\ \hline
\end{tabular}
\caption{ Error analysis \eqref{NonlinearScalar}.}
\end{table}

%###############################################
%\newpage
\subsection{Nonlinear systems of equations}

Let us consider the following system of nonlinear
integral equations
\begin{equation}\label{NonlinearSystem-01}
 \begin{cases}
    \int\limits_{0}^{\frac{t}{2}} s(t+s) x^2_1(s))ds + \int\limits_{\frac{t}{2}}^{t}(3x_2(s)+x^3_2(s))ds= f_1(t),\\
    \int\limits_{0}^{\frac{t}{2}} (1+t-s)( x_1(s)-x^2_1(s))ds - \int\limits_{\frac{t}{2}}^{t}( x_2(s) + x^4_2(s))ds= f_2(t),
\end{cases}
\end{equation}
where \(t\in [0,1]\), and functions \(f_i(t), \; i=1,2,\) are selected such as
\[
  x_1(t) = t^2,\; x_2(t) =t^3
\]
is solution of the  system.
Let us select the initial approximation far enough from the exact one
\[   x^0_1(t) = 0.4 t^2,\; x^0_2(t) =0.5 t^3. \]
Results are given in Tab. 4, where \(N_{it}\) is iterations number, \(m\)
is order of approximating polynomials.
\begin{table}[h!]\label{Num-4}
 \centering
\begin{tabularx}{\textwidth}{|>{\hsize=0.075\textwidth}>{\centering}X|>{\hsize=0.075\textwidth}>{\centering}X|>{\centering}X|>{\centering}X|>{\centering}X|}
\hline
    \(N_{it}\)  & \(m\)  & \(\varepsilon_1\)      & \(\varepsilon_2\)     & \(\varepsilon\)\tabularnewline
\hline
     1     & 3      &\(0.426435\)             &\(0.133873\)             &\(0.446955\) \tabularnewline
\hline
     10    & 3      &\(5.97261\cdot 10^{-5}\) &\(1.91177\cdot 10^{-5}\) &\(6.27112\cdot 10^{-5}\) \tabularnewline
\hline
     20    & 3      &\(2.81636\cdot 10^{-9}\) &\(9.78086\cdot 10^{-10}\)&\(2.98137\cdot 10^{-9}\) \tabularnewline
\hline
\end{tabularx}
 \caption{Results for problem \eqref{NonlinearSystem-01}.}
\end{table}

$\,$\\

Let us now consider the same system, but not exponential exact solution
\begin{equation}\label{NonlinearSystem-02}
 \begin{cases}
    \int\limits_{0}^{\frac{t}{2}} s(t+s) x^2_1(s)ds + \int\limits_{\frac{t}{2}}^{t}(3 x_2(s)+x^3_2(s))ds= f_1(t),\\
    \int\limits_{0}^{\frac{t}{2}} (1+t-s)( x_1(s)-x^2_1(s))ds - \int\limits_{\frac{t}{2}}^{t}( x_2(s) + x^4_2(s))ds= f_2(t),
\end{cases}
\end{equation}
where \(t\in [0,1]\) and functions \(f_i(t), \; i=1,2,\) are selected such as the following functions are exact solutions of the example under consideration
\[
  x_1(t) = \cos(t), \; x_2(t) = \sin(t).
\]

Results are given in Tab. 5 and 6.

\begin{table}[h!]\label{Num-5}
 \centering
\begin{tabularx}{\textwidth}{|>{\hsize=0.075\textwidth}>{\centering}X|>{\hsize=0.075\textwidth}>{\centering}X|>{\centering}X|>{\centering}X|>{\centering}X|}
\hline
    \(N_{it}\)  & \(m\)  & \(\varepsilon_1\)      & \(\varepsilon_2\)     & \(\varepsilon\)\tabularnewline
\hline
     1     & 5      &\( 0.156589\)             &\( 0.0170903\)             &\(0.157519\) \tabularnewline
\hline
     5     & 5      &\(5.04072\cdot 10^{-4}\)  &\(4.89738\cdot 10^{-5}\)   &\(5.06445\cdot 10^{-4}\) \tabularnewline
\hline
     10    & 5      &\(5.22785\cdot 10^{-6}\) &\(4.16019\cdot 10^{-6}\) &\(6.68114\cdot 10^{-6}\) \tabularnewline
\hline
     20    & 5      &\(1.97551\cdot 10^{-6}\) &\(3.69009\cdot 10^{-6}\)&\(4.18562\cdot 10^{-6}\) \tabularnewline
\hline
\end{tabularx}
 \caption{Results for problem \eqref{NonlinearSystem-02} for \(m=5\).}
\end{table}
%*********************************
\begin{table}[h!]\label{Num-6}
 \centering
\begin{tabularx}{\textwidth}{|>{\hsize=0.075\textwidth}>{\centering}X|>{\hsize=0.075\textwidth}>{\centering}X|>{\centering}X|>{\centering}X|>{\centering}X|}
\hline
    \(N_{it}\)  & \(m\)  & \(\varepsilon_1\)      & \(\varepsilon_2\)     & \(\varepsilon\)\tabularnewline
\hline
     1     & 10      &\( 0.194802\)             &\(0.0224296\)             &\(0.196089\) \tabularnewline
\hline
     5     & 10      &\(3.02512\cdot 10^{-3}\)             &\(3.21621\cdot 10^{-4}\)             &\(3.04217\cdot 10^{-3}\) \tabularnewline
\hline
     10    & 10      &\( 2.28126\cdot 10^{-5}\) &\(2.22376\cdot 10^{-6}\) &\(2.29207\cdot 10^{-5}\) \tabularnewline
\hline
     20    & 10      &\(4.03156\cdot 10^{-9}\) &\(6.04141\cdot 10^{-10}\)&\(4.07657\cdot 10^{-9}\) \tabularnewline
\hline
\end{tabularx}
 \caption{Results for problem  \eqref{NonlinearSystem-02} for \(m=10\).}
\end{table}

Here, to improve the convergence of the method, the initial approximation was chosen sufficiently close to the exact one:
\[
  x^0_1(t) = 0.9 \cos(t),\; x^0_2(t) = 0.9 \sin(t).
\]

%---------------------------------------------------------------------
%----------------------  S E C T I O N 5 -----------------------------
%---------------------------------------------------------------------
\section{Conclusion}

In this article the efficient  numerical methods for solution of the novel class of weakly regular linear and nonlinear systems of Volterra integral
equations of the first kind are proposed.
The iterative method for solution of the system of nonlinear Volterra integral
equations of the first kind with special jump discontinuous kernels employs the  linearization of the integral vector-operator using the  Newton-Kantorovich iterative process.
The convergence theorem for proposed iterative process is formulated. Polynomial collocation method
is also proposed to solve such systems.
The efficiency constructed theory and numerical methods is illustrated in five examples.

\section*{Acknowledgment}

This work is carried out partly within the framework of the research projects III.17.3.1, III.17.4 of the program of fundamental research of the Siberian Branch of the Russian Academy of Sciences, reg No. AAAA-A17-117030310442-8.



%% References


\bibliographystyle{elsarticle-num}

\begin{thebibliography}{00}



\harvarditem{Kantorovich \& Gorkov}{1959}{kantor}
\textsc{L. Kantorovich, L. Gorkov, }
 On some functional equations arising in analysis of single-commodity economic model, {\em Dokl Akad. Nauk SSSR,} \textbf{129} (4) (1959) 732--736.


\harvarditem{Hritonenko \& Yatsenko}{1996}{hriton}
\textsc{N. Hritonenko, Yu.  Yatsenko, }   {\em Modeling and Optimization of the Lifetime of Technologies.}   Dordrecht: Kluwer Academic Publishers, 1996.

\harvarditem{Apartsyn}{2003}{apartsyn}
\textsc{A.S. Apartsyn,}   {\em Nonclassical Linear Volterra Equations of the First Kind.}  De Gruyter, 2003.



\harvarditem{Sidorov}{2014}{Sidorov-book}
\textsc{D. Sidorov,}  {Integral Dynamical Models: Singularities, Signals and Control.}
{In:\textsc{L. O. Chua,} ed. \em  World Scientific Series on Nonlinear Sciences Series A:} Vol. 87, Singapore: World Scientific Press, 2015, 243~p.


\harvarditem{Muftahov, Sidorov et al.}{2016}{arxiv}
\textsc{I. Muftahov, D. Sidorov, A. Zhukov, D. Panasetsky, A. Foley, Y. Li, A. Tynda,}  {Application of Volterra Equations to Solve Unit Commitment Problem of Optimised Energy Storage and Generation.}
arXiv:1608.05221.

\harvarditem{Kythe \& Puri}{2002}{book}
\textsc{P. K. Kythe, P. Puri, }   {\em Computational Methods for Linear Integral Equations.}   Birkh\"{a}user. Boston, 2002.

\harvarditem{Evans}{1910}{Evans}
\textsc{G. C. Evans,} {
Integral {e}quation of the {s}econd {k}ind with discontinuous {k}ernel,} Transactions of the American Mathematical
Society. \textbf{11}
(4) (1910)  393--413.



\harvarditem{Anselone}{1967}{Anselone}
\textsc{P. M. Anselone,} {
 {U}niform approximation {t}heory for integral {e}quations with
{d}iscontinuous {k}ernels,} SIAM J. Numer. Anal. \textbf{4} (1967) 245--253.

  \harvarditem{Denisov \& Lorenzi}{1995}{denisovlorenzi}
\textsc{A. M. Denisov, A. Lorenzi.}  On a special {V}olterra integral equation of the first kind. {\em {B}ollettino {U}.{M}.{I}.}, \textbf{7} (9-B) (1995) 443--457.

\harvarditem{Brunner}{1997}{brunner}
\textsc{H. Brunner,}
 1896 {--} 1996: One hundred years of Volterra integral equations
of the first kind, {\em Applied Numerical Mathematics,} \textbf{24} (1997)
83--93.


\harvarditem{Sidorov}{2011}{Sidorov-FR}
\textsc{D. N. Sidorov, }  Volterra  equations of the first kind with
discontinuous kernels in the
theory of evolving systems control. {\em Stud. Inform. Univ.}, \textbf{9} (2011) 135--146.




\harvarditem{Sidorov}{2013}{SidorovDU}
\textsc{D. N. Sidorov, }
On parametric families of solutions of Volterra integral equations of the first kind with piecewise smooth kernel, {\em Differential Equations,} \textbf{49}(2) (2013) 210--216.


\harvarditem{Sidorov}{2013}{Sidorov-2013}
\textsc{D. N. Sidorov, } Solvability of Systems of Volterra Integral Equations of the First Kind with Piecewise Continuous Kernels,
 {\em Russian Mathematics (Iz. VUZ)}, \textbf{57}(1) (2013) 54--63.
%------------------


\harvarditem{Markova \& Sidorov}{2014}{SidorovMarkovaAiT}
\textsc{E. V. Markova, D. N. Sidorov, }
On one integral Volterra model of developing dynamical systems,
{\em Automation and Remote Control,}
 \textbf{75} (3) (2014)  413--421.


\harvarditem{Sidorov \& Tynda \& Muftahov}{2014}{vestnik}
\textsc{D. N. Sidorov, A. N. Tynda, I. R. Muftahov  }
 Numerical solution of the Volterra integral equations of the first kind with piecewise continuous kernel, {\em Bulletin of the South Ural State University, Series: Mathematical Modelling, Programming and Computer Software}, \textbf{7} (3) (2014) 107--115.


 \harvarditem{Lorenzi}{2013}{Lorenzi2013}
\textsc{A. Lorenzi, }  {O}perator equations of the first kind and integro-differential equations of degenerate type in {B}anach spaces and applications to integro-differential {P}{D}{E}s. {\em Eurasian Journal of Mathematical and Computer Applications}, \textbf{1}(2) (2013) 50--75.


   \harvarditem{Sidorov \& Sidorov}{2014}{sidsid}
\textsc{N.A. Sidorov, D.N. Sidorov,}  On the solvability of a class of Volterra operator equations of the first kind with piecewise continuous kernels. {\em Math. Notes}, \textbf{96}(5) (2014) 811--826.



\harvarditem{Khromov}{2006}{Khromov}
\textsc{A. P. Khromov,} {
Integral {o}perators with {d}iscontinuous {k}ernel on {p}iecewise {l}inear {c}urves,} Sbornik: Mathematics. \textbf{197}
(11) (2006) 115--142.


\harvarditem{Davies \& Duncan}{2017}{Davies}
\textsc{ P.J. Davies, D.B. Duncan}
Numerical approximation of first kind Volterra convolution integral equations with discontinuous kernels, {\em J. Integral Equations Appl.}, \textbf{28} (1) (2017) 41--73.



  \harvarditem{Boikov \& Tynda}{2003}{Boikov-Tynda11}
  \textsc{I.V. Boikov, A. N. Tynda, }  Approximate solution of nonlinear integral equations of
  developing systems theory. {\em Differential Equations}, \textbf{39} (9) (2003) 1214--1223.

  \harvarditem{Tynda}{2006}{Tynda-18}
   \textsc{A. N. Tynda, }   Numerical  algorithms of optimal complexity for weakly singular
     Volterra integral equations, {\em Comp. Meth. Appl. Math.},
     \textbf{6} (4) (2006) 436--442.

\harvarditem{Muftahov \& Sidorov}{2016}{Muftahov-Ural}
\textsc{I.R. Muftahov,  D.N. Sidorov, }  Solvability and numerical solutions of systems of nonlinear Volterra integral equations of the first kind with piecewise continuous kernels. {\em Vestnik YuUrGU. Ser. Mat. Model. Progr.}, \textbf{9}(1) (2016) 130--136.

	


\harvarditem{Muftahov et al}{2017}{Muftahov-CAM}
\textsc{I. Muftahov, A. Tynda, D. Sidorov, }  Numeric solution of Volterra integral equations of the first kind with discontinuous kernels. {\em Journal of Computational and Applied Mathematics}, \textbf{313} (2017) 119--128














\harvarditem{Kantorovich \& Akilov}{1982}{Kantorovich2}
  \textsc{L. V. Kantorovich, G. P. Akilov,}  {\em Functional Analysis}. Pergamon; 2nd edition, 1982, 589~p.














%\harvarditem{Press}{2007}{Recipes}
%\textsc{W. H. Press, S. A. Teukolsky, W. T. Vetterling, B. P. Flannery, }  Numerical recipes in C: the art of scientific computing.
%{ \em Cambridge University Press}, 2nd ed. 2007.
%
%\harvarditem{Solow}{1969}{solow}
%\textsc{R. M. Solow, }   {Investment and technical progress.}  In: \textsc{ K. J. Arrow, S. Karlin and P. Suppes,} ed.
% {\em Mathematical Methods in the Social Sciences,} Stanford University
%Press, 1969,  pp.~89--104.


























%------------------



%------------------

%------------------
%  \bibitem{Tynda-pamm1}
%   A.N. Tynda,  Numerical methods for 2D weakly singular Volterra integral equations
%                of the second kind. \emph{PAMM}, Volume 7 (2007), Issue 1.
%------------------








%\harvarditem{Arnold \textit{et al.}}{1997}{arnoldthree}
%\textsc{Arnold, D.N., Falk, R.S. \& Winther, R.} (1997)   Preconditioning
%in $H(div)$ and applications,   {\em Math. Comp.,}   \textbf{66}
%(219), 957--984.
%
%\harvarditem{Arnold \textit{et al.}}{2000}{arnoldtwo}
%\textsc{Arnold, D.N., Falk, R.S. \& Winther, R.} (2000)   Multigrid
%in $H(div)$ and $H(curl)$,   {\em Numer. Math.,}   \textbf{85},
%197--218.
%
%\harvarditem{Brezzi}{1974}{brezzi1}
%\textsc{Brezzi, F.} (1974)  {On the existence, uniqueness and approximation
%of saddle point problems arising from Lagrangian multipliers,}
% {\em RAIRO Anal. Num\'er.,} \textbf{8}, 129--151.



%
%\harvarditem{Brezzi \& Bathe}{1990}{brezzibathe}
%\textsc{Brezzi, F. \& Bathe, K.J.} (1990)  {A discourse on the stability
%conditions for mixed finite element formulations},  {\em Comp. Methods.
%Appl. Mech. Engrg.,}   \textbf{82}, 27--57.
%
%\harvarditem{Renardy \& Rogers}{2004}{bit}
%\textsc{Renardy, M. \& Rogers, R.C.} (2004)  {\em An introduction to partial differential % equations.}   New York: Springer-Verlag.


%
%\harvarditem{Cai \textit{et al.}}{1993}{cai2}
%\textsc{Cai, Z., Goldstein, C.I. \& Pasciak, J.E.} (1993)  {Multilevel
%iteration for mixed finite element systems with penalty},  {\em SIAM
%J. Sci. Comput.,}   \textbf{14} (5), 1072--1088.
%






%
%\harvarditem{Ewing \& Wang}{1992}{ewingone}
%\textsc{Ewing, R.E. \& Wang, J.} (1992)  {Analysis of the Schwarz
%algorithm for mixed finite elements methods},  {\em M$^2$AN, Math.
%Model. Numer. Anal.,}   \textbf{26} (6), 739--756.
%
%\harvarditem{Ewing \& Wang}{1994} {ewingtwo}
%\textsc{Ewing, R.E. \& Wang, J.} (1994)  {Analysis of multilevel
%decomposition iterative methods for mixed finite element methods},  {\em M$^2$AN, Math. Model. Numer. Anal.,}   \textbf{28} (4), 377--398.
%
%
%\harvarditem{Fortin \& Glowinski}{1983}{glow}
%\textsc{Fortin, M. \& Glowinski, R.} (1983)  {\em Augmented Lagrangian
%Methods: Applications to the Numerical Solution of Boundary-Value Problems.}
%Amsterdam: North-Holland.
%
%\harvarditem{Hiptmair}{1997}{hiptmair}
%\textsc{Hiptmair, R.} (1997)  {Multigrid methods for $H(div)$ in
%three dimensions},  { \em Electron. Trans. Numer. Anal.,}
%\textbf{6}, 133--152.
%
%\harvarditem{Kellogg}{1975}{kellogg}
%\textsc{Kellogg, R.B.,} (1975)  {On the Poisson equation with intersecting
%interfaces},  {\em Appl., Anal.,}   \textbf{4}, 101--129.
%
%\harvarditem{Morin \textit{et al.}}{2002}{morin}
%\textsc{Morin, P., Nochetto, R.H. \& Siebert, K.G.} (2002)  {Convergence
%of adaptive finite element methods},  {\em SIAM Rev.,}
%1--28.
%
%\harvarditem{Paige \& Saunders}{1975}{paige}
%\textsc{Paige, C.C. \& Saunders, M.A.} (1975)  {Solution of sparse
%indefinite systems of linear equations,}  {\em SIAM J. Numer. Anal.,}
%  \textbf{12}, 617--629.
%
%\harvarditem{Powell \& Silvester}{2003}{me1}
%\textsc{Powell, C.E. \& Silvester, D.} (2003)  {Optimal preconditioning
%for Raviart-Thomas mixed formulation of second-order elliptic problems,}
% {\em SIAM J. Matrix Anal. Appl.,} \textbf{25} (3), 718--738.
%
%\harvarditem{Raviart \& Thomas}{1977}{ravtom}
%\textsc{Raviart, P.A. \& Thomas, J.M.} (1977)  {A mixed finite element
%method for second order elliptic problems.}  {In: \textsc{I. Galligani
%and E. Magenes,} eds. \em Mathematical Aspects of the Finite Element Method,
%Lect. Notes in Math.,} 606, New York: Springer-Verlag, 292--315.
%
%\harvarditem{Russell \& Wheeler}{1983}{russell}
%\textsc{Russell, T.F. \& Wheeler, M.F.} (1983)  {Finite element and
%finite difference methods for continuous flows in porous media.}  {In:
%\textsc{R.E. Ewing,} ed. \em  The Mathematics of Reservoir Simulation.} Philadelphia:
%SIAM, 35--106.
%
%\harvarditem{Rusten \& Winther}{1992}{rusten}
%\textsc{Rusten, T. \& Winther, R.} (1992)  {A preconditioned iterative
%method for saddlepoint problems},  {\em SIAM J. Matrix. Anal. Appl.,}
%  \textbf{13} (3), 887--904.
%
%\harvarditem{Rusten \textit{et al.}}{1996}{rusten3}
%\textsc{Rusten, T., Vassilevski, P.S. \& Winther, R.} (1996)  {Interior
%preconditioners for mixed finite element approximations of elliptic problems},
% {\em Math. Comp.,}   \textbf{65} (214), 447--466.
%
%\harvarditem{Scheichl}{2000}{scheichl}
%\textsc{Scheichl, R.} (2000)  {\em Iterative solution of saddle point
%problems using divergence-free finite elements with applications to groundwater
%flow.} Thesis (PhD). Bath University.
%



%------------------



%
%\harvarditem{Vassilevski \& Wang}{1992}{vasswang}
%\textsc{Vassilevski, P.S. \& Wang, J.} (1992)  {Multilevel iterative
%methods for mixed finite element discretizations of elliptic problems},  {\em
%Numer. Math.},   \textbf{63}, 503--520.
%
%\harvarditem{Vassilevski \& Lazarov}{1996}{vasslaz}
%\textsc{Vassilevski, P.S. \& Lazarov, R.D.} (1996)  {Preconditioning
%mixed finite element saddle point elliptic problems,}  {\em Numer.
%Linear Algebra Appl.,}   \textbf{3} (1), 1--20.
\end{thebibliography}


 
%%%%%%%%%%%%%%%%%%%%%%%%

%\appendix
%\begin{Appendix}
%\section{}
%\section*{Appendix}
%In this section we want to show the tables which consist of the results obtained by our methods.

\end{document}